\documentclass[a4paper]{article}
\usepackage[margin=1.5in]{geometry}
\usepackage[english]{babel}
\usepackage[utf8]{inputenc}
\usepackage{tabularx}
\usepackage[T1]{fontenc}
\usepackage[utf8]{inputenc}
\usepackage{amssymb}
\usepackage{graphicx}
\usepackage[colorinlistoftodos]{todonotes}
\usepackage{tikz}
\usepackage{authblk}
\usepackage{amsmath}

\author[1]{U S Naveen Balaji}
\author[2]{S Sivasankar}
\author[3]{Sujan Kumar S}
\author[4]{Vignesh Tamilmani}
\affil[1,3,4]{Department of Science and Humanities, PES University, Bangalore - 560085, Karnataka, India. email(s):  naveenandmetallica@gmail.com, sivshankar@gmail.com, sujankumar175@gmail.com, vignesht16202@gmail.com}
\affil[2]{Department of Mathematics, R V Institute of Technology and Management, Bangalore - 560059, Karnataka, India.}
\title{Cyclic Symmetry of Riemann Tensor in Fuzzy Graph Theory}

\date{\vspace{-5ex}}
\begin{document}
\maketitle

\begin{abstract}
\noindent
In this paper, we define a graph-theoretic analog for the Riemann tensor and analyze properties of the cyclic symmetry. We have developed a fuzzy graph theoretic analog of the Riemann tensor and have analyzed its properties.  We have also shown how the fuzzy analog satisfies the properties of the $6\times 6$ matrix of the Riemann tensor by expressing it as a union of the fuzzy complete graph formed by the permuting vertex set and a Levi-Civita graph analog. We have concluded the paper with a brief discussion on the similarities between the properties of the fuzzy graphical analog and the Riemann tensor and how it can be a plausible analogous model for the Petrov-Penrose classification.
\end{abstract}

\section{Introduction}
Tensors and Differential geometry are central to General Relativity, they are the foundation to the seminal theory of Einstein. The Riemann curvature tensor named after Bernhard Riemann is a higher-dimensional analogue of the Gaussian
curvature and is closely related to tidal forces, it represents
the tidal force experienced by a particle moving along a geodesic. In $4$-dimensions, the Riemann tensor has $256$ components and observations reveal a variety of algebraic symmetries such as the first skew symmetry, the second skew symmetry, and the block symmetry all of which reduce the $256$ components to $20$ independent components. The last algebraic symmetry, called the cyclic symmetry is closely associated with Bianchi's first identity. For values of $n>3$ in $\frac{1}{2}n^{2}(n^{2}-1)$ algebraically independent components of the Riemann tensor the components are represented by the Weyl tensor which also possesses all three algebraic symmetries and in addition it can be thought of as that part of the curvature tensor such that all contractions vanish, i.e., a pseudo-Riemannian manifold is said to be conformally flat if its Weyl tensor vanishes. It was A.Z. Petrov who classified the algebraic symmetries of the Weyl tensor, called the Petrov-Penrose classification. Generally, gravitational fields are classified in accordance to the Petrov-Penrose classification of their corresponding Weyl tensor.
\\\\ 
In this paper we develop a graph theoretic analog of the Riemann tensor which we then use to help develop a fuzzy graph analog of the Petrov-Penrose classification. We exploit the cyclic symmetry of the Riemann tensor to help define the graphical analog and also discuss, through theorems, the similarities and properties of the graphical analog to its tensor form. The paper is organized as follows. Section $2$ contains the preliminaries and in section $3$ fuzzy graphical analogs of the Riemann tensor and the Levi-Civita symbol are defined. Section $4$ deals with the fuzzy approach to Pentrov-Penrose classification.  

\section{Preliminaries}
\label{sec:introduction}
\textbf{Definition1.1}[1]. Let $R^{n}$ denote the Euclidean space of $n$-dimensions, i.e., the set of all $n$-tuples $\left(x^{1},x^{2},...,x^{n} \right)$ $\left(-\infty < x^{i} < \infty  \right)$ with the usual topology (open and closed sets are defined in the usual way), and let $\frac{{R}^{n}}{2}$ denote the \textit{lower half} of ${R}^{n}$, i.e. the region of ${R}^{n}$ for which $x^{1} \leq 0$. A map $\phi$ of an open set $\mathcal{O} \subset {R}^{n}$ (respectively $\frac{{R}^{n}}{2}$) to an open set $\mathcal{O}' \subset {R}^{n}$ (respectively $\frac{{R}^{n}}{2}$) is said to be of class $C^{r}$ if the coordinates $\left(x'^{1},x'^{2},...,x'^{m} \right)$ of the image point $\phi(p)$ in $\mathcal{O}'$ are $r$-times continuously differentiable
functions of the coordinates $\left(x'^{1},x'^{2},...,x'^{m} \right)$ of $p$ in $\mathcal{O}$. 
\\

\noindent
\textbf{Definition 1.2}[1]. A $C^{r}$ $n$-dimensional manifold $\mathcal{M}$ is a set $\mathcal{M}$ together with a $C^{r}$ atlas $\{\mathcal{U}_{\alpha},\mathcal{\phi}_{\alpha} \}$, where the $\mathcal{U}_{\alpha}$ are subsets of $\mathcal{M}$ and the $\phi_{\alpha}$ are one-one maps of the corresponding $\mathcal{U}_{\alpha}$ to open sets in ${R}^{n}$ such that\\
\newline
\noindent
$(1)$ The $\mathcal{U}_{\alpha}$ cover $\mathcal{M}$, i.e. $\mathcal{M}= \bigcup_{\alpha}\mathcal{U}_{\alpha}$, \\

\noindent
$(2)$ if $\mathcal{U}_{\alpha}\cup \mathcal{U}_{\beta}$ is non-empty, then the map $\phi_{\alpha} \circ \phi_{\beta}^{-1}: \phi_{\beta} \left(\mathcal{U}_{\alpha} \cup \mathcal{U}_{\beta} \right) \rightarrow \phi_{\alpha}\left(\mathcal{U}_{\alpha} \cup \mathcal{U}_{\beta} \right)$ is a $C^{r}$ map of an open subset of ${R}^{n}$ to an open subset of ${R}^{n}$. \\

\noindent
\textbf{Definition 1.3}[2,7]. The mapping $f:V\rightarrow U$, where the open sets $V,U \in{{R}^{n}}$, is called a homeomorphism if it is bijective and if $f$ and its inverse $f^{-1}$ are continuous. \\

\noindent
\textbf{Definition 1.4}[2,7]. A chart for a topological space $\mathcal{M}$ is a homeomorphism $\phi$ from an open subset $U$ of $\mathcal{M}$ to an open subset of a Euclidean space. The chart is traditionally recorded as the ordered pair $\left(U, \phi \right)$.\\

\noindent
\textbf{Definition 1.5}[2,7]. A tangent vector $v$ to the differential manifold $\mathcal{M}$ at a point $p \in \mathcal{M}$ is defined as $(\left(V_{\rho}, z_{\rho}, v_{z_{\rho}} \right)$, where $(\left(V_{\rho}, z_{\rho}\right)$ are charts which contain $p$ and $v_{z_{\rho}} = v_{z_{\rho}}^{j}$, $j = 1,2,...,n$ are vectors in ${R}^{n}$.\\

\noindent
\textbf{Definition 1.6}[2,7]. Let $\Lambda$ be a $p$-form field defined by $\Lambda = \Lambda_{\alpha \beta... \zeta}dx^{\alpha} \wedge dx^{\beta}\wedge ... \wedge dx^{\zeta}$ where $\alpha, \beta, ..., \zeta$ are arbitrary indices. The exterior derivative acts on this $p$-form field to produce a $(p+1)$-form field as follows
\\
\begin{equation}
    d\Lambda = d\Lambda_{\alpha \beta... \zeta}dx^{\alpha} \wedge dx^{\beta}\wedge ... \wedge dx^{\zeta}.
\end{equation}
\\
The exterior derivative of a $(p-1)$-form produces a $p$-form defined as follows
\\
\begin{equation}
    \Xi = \frac{1}{p!}\Xi_{i_{1}\ i_{2}\ ...\ i_{p}}dx^{i_{1}}\wedge dx^{i_{2}}\wedge ... \wedge dx^{i_{p}}
\end{equation}
\\
and the exterior derivative is defined as
\\
\begin{equation}
    d\Xi = \frac{1}{p!}\frac{\partial \Xi_{i_{1}\ i_{2}\ ...\ i_{p}}}{\partial x^{i_{0}}}dx^{i_{0}}\wedge dx^{i_{1    }}\wedge ... \wedge dx^{i_{p}}.
\end{equation}
\\

\noindent
\textbf{Definition 1.7}[1]. A Cartesian product is defined as the ordered set of vectors and one-forms $\left(\eta^{1},...,\eta^{m},\textbf{Y}_{1},...,\textbf{Y}_{n} \right)$, where the $\textbf{Y}$'s and $\eta$'s are arbitrary vectors and one-forms respectively. The Cartesian product is expressed as the product of the tangent space $T_{p}$ of vectors at a point $p$ and the tangent space's dual or the cotangent space $^{*}T$ of $1$-forms at $p$ written as follows
\\
\begin{equation}
    \Pi_{m}^{n} = \underbrace{^{*}T_{p}\times^{*}T_{p}\times...\times^{*}T_{p}\times}_{n\ factors}\times \underbrace{T_{p}\times T_{p}\times...\times T_{p}}_{m\ factors}.
\end{equation}
\\

\noindent
\textbf{Definition 1.8}[1]. A tensor of rank $\left(\begin{array}{lr} n \\ m \end{array}\right)$ at a point $p$ is a function on $\Pi_{m}^{n}$ which is linear in each argument, i.e., if $\textbf{T}$ is a tensor of rank $\left(\begin{array}{lr} n \\ m \end{array}\right)$at $p$, the number into which $\textbf{T}$ maps the element $\left(\eta^{1},...,\eta^{m},\textbf{Y}_{1},...,\textbf{Y}_{n}\right)$ of $\Pi_{m}^{n}$ as $T\left(\eta^{1},...,\eta^{m},\textbf{Y}_{1},...,\textbf{Y}_{n}\right)$, where the $\textbf{Y}$'s and $\eta$'s are arbitrary vectors and one-forms respectively.\\

\noindent
\textbf{Definition 1.9}[2,3]. The Riemann tensor is a four-index tensor which has $256$ components in $4$-dimensions. Making use of the symmetry relations,
\\
\begin{equation}
    R_{iklm} = -R_{ikml} =- R_{kilm},
\end{equation}
\\
the number of independent components is reduced to $36$. Using the condition
\\
\begin{equation}
    R_{iklm} = R_{lmik},
\end{equation}
\\
the number of coordinates reduces to $21$. Finally, using 
\\
\begin{equation}
    R_{iklm} + R_{ilmk} + R_{imkl} = 0,
\end{equation}
\\
$20$ independent components are left.
\\

\noindent
\textbf{Definition 1.10}. A graph $G = (V,E)$ is analogous to the Riemann tensor with a vertex set $V = \{v_{1},v_{2},v_{3},v_{4} \}$, where $v_{1} = i$, $v_{2} = k$, $v_{3} = l$, $v_{4} = m$, if it satisfies the following:
\\\\
\noindent
a. The vertex $v_{1}$ is fixed and is connected to only one in $V-\{v_{1}\}$, \\
b. There are three vertices which span a $K_{3}$ graph amongst themselves, \\
c. The direction of the cycle in the $K_{3}$ graph determines the overall sign assigned to the graph, i.e., if the cycle is in the counter clockwise direction, we assign a positive sign and if the cycle is in the clockwise direction, we assign a negative sign.\\

\begin{figure}\label{fig:Fuzzy1}
    \centering
    \includegraphics{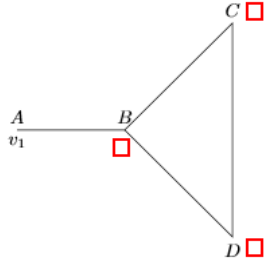}
    \caption{The vertex points that occupy the vertex positions $B$, $C$, and $D$ for the subgraph $H$ $(=K_{3})$, where $H\subset{G}$ and the occupation of $v_{2}$, $v_{3}$, and $v_{4}$ in each of these positions is equally likely.}
    \label{fig:my_label}
\end{figure}

\begin{figure}
    \centering
    \includegraphics[height=10cm, width=8.5cm]{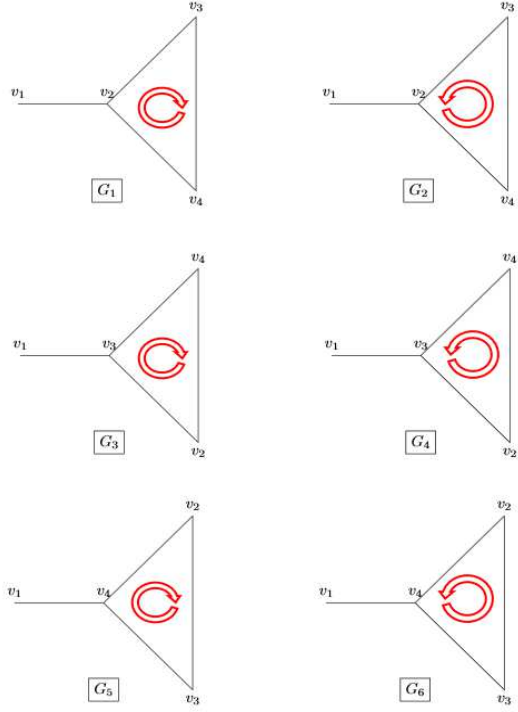}
    \caption{From the above figure, $G_{1}=-G_{2}$, $G_{3} = -G_{4}$, and $G_{5} = -G_{6}$}
    \label{fig:Fuzzy2}
\end{figure}

\noindent
\textbf{Definition 1.11}. Let $H \subseteq G$ denote the $K_{3}$ graph in which the three positions are labelled $B, C$ and $D$, with $B$ being adjacent to the fixed vertex position $A$. We note that in position $B$, the vertices $v_{2}, v_{3}$ and $v_{4}$ are equally likely for occupation and once a particular vertex is occupied, the others vertices occupy positions $C$ and $D$ in cyclic order (See Figure \ref{fig:Fuzzy1}).
\\
\\
Suppose $v_{2}$ occupies the vertex position $B$ then, $v_{3}$ would occupy $C$ and $v_{4}$ would occupy $D$. For each arrangement of the vertices we have two variants which differ from each other by a negative sign which is determined by the direction of the cycle assigned to the subgraph $H$. Thus, for each combination of the adjacent vertices we have three variants listed below (see Figure \ref{fig:Fuzzy2}):\\\\
\noindent
1. With $\{v_{1},v_{2},v_{3},v_{4}\}$ we have graph $G_{1}$ for clock wise cycle and $G_{2}$ for a counter clock wise cycle which are related as, $G_{1} = -G_{2}$.\\
2. With $\{v_{1},v_{3},v_{4},v_{2}\}$ we have graph $G_{3}$ for clock wise cycle and $G_{4}$ for a counter clock wise cycle which are related as, $G_{3} = -G_{4}$.\\
3. With $\{v_{1},v_{4},v_{2},v_{3}\}$ we have graph $G_{5}$ for clock wise cycle and $G_{6}$ for a counter clock wise cycle which are related as, $G_{5} = -G_{6}$.\\
\\
\\
\noindent
\textbf{Definition 1.12}. Let $T(\alpha,\beta,\gamma,\delta) = T_{\alpha \beta \gamma \delta}$ denote a function which takes the sequence of the vertices as an input and displays the vertices as indices of a particular graph $T$ as the output. With this formulation, we can express a graph $G_{i}$ as follows: \\
\noindent
a. $G_{1} \equiv G\left(v_{1}, v_{2}, v_{3}, v_{4} \right) = G_{v_{1}v_{2}v_{3}v_{4}} = G_{iklm}$, \\
b. $G_{2} \equiv G\left(v_{1}, v_{2}, v_{4}, v_{3} \right) = G_{v_{1}v_{2}v_{4}v_{3}} = G_{ikml}$, \\
c. $G_{3} \equiv G\left(v_{1}, v_{3}, v_{4}, v_{2} \right) = G_{v_{1}v_{3}v_{4}v_{2}} = G_{ilmk}$, \\
d. $G_{4} \equiv G\left(v_{1}, v_{3}, v_{2}, v_{4} \right) = G_{v_{1}v_{3}v_{2}v_{4}} = G_{ilkm}$, \\
e. $G_{5} \equiv G\left(v_{1}, v_{4}, v_{2}, v_{3} \right) = G_{v_{1}v_{4}v_{2}v_{3}} = G_{imkl}$, \&\\
f. $G_{6} \equiv G\left(v_{1}, v_{4}, v_{3}, v_{2} \right) = G_{v_{1}v_{4}v_{3}v_{2}} = G_{imlk}$. \\
\noindent
\textbf{Definition 1.13}. Let $G_{ik(lm)} = \frac{1}{2!}\left(G_{iklm} - G_{ikml} \right)$ such that 
\\
\begin{equation}
  2\left[G_{ik(lm)} + G_{il(mk)} + G_{im(kl)} \right] = 0  
\end{equation}
\\

\noindent
\textbf{Theorem 1.1}. Three permuting indices, by definition $1.10$, condenses six graphs to three graphs and thus, 
\\
\begin{equation}
    G_{i(klm)} = 0.
\end{equation}
\\
\textit{Proof} To condense the number of graphs we can make use of the fact that $k, l$ and $m$ permute cyclically as follows
\\
\begin{equation}
    G_{i(klm)} = \frac{1}{3!}\left[G_{iklm} + G_{ilmk} + G_{imkl} \right],
\end{equation}
\\
and thus, $3!\left[G_{iklm} + G_{ilmk} + G_{imkl} \right] = 0$ or $ G_{i(klm)} = 0$. \\

\noindent
\textbf{Theorem 1.2}. The antisymmetry in each pair of indices (vertices) of a graph $G$ constructed from definition $1.10$ implies that there are $P = \frac{1}{2}n\left(n-1\right)$ ways of choosing independent pairs of indices. \\
\noindent
\textit{Proof} In our discussion we have $n=4$ indices and thus there are $P = \frac{1}{2}(4)(3) = 6$ ways of choosing pairs which resulted in graphs $G_{1}, G_{2}, G_{3},G_{4}, G_{5},$ and $G_{6}$.\\

\noindent
\textbf{Theorem 1.3}. A graph $G$ with one fixed vertex and $\alpha$ number of permuting vertices possesses $1 + \xi(K_{\alpha})$ number of edges, where $\xi(K_{\alpha})$ is the number of edges of a $K_{\alpha}$ complete graph.\\
\noindent
\textit{Proof} Since every graph has $\xi \leq \frac{\gamma(\gamma - 1)}{2}$ number of edges where $\gamma$ is the number of vertices, the number of edges of a graph constructed from definition $1.10$ will be the number of edges of the subgraph $H$ and the edge connecting the fixed vertex to one of the vertices of the permuting vertex set. Since the subgraph $H$ is a $K_{3}$ graph, we have $1 + \xi(K_{3})$ number of edges.\\

\noindent
\textbf{Definition 1.14}. Let $u_{1} = ik, u_{2} = il, u_{3} = im, u_{4} = kl, u_{5} = km, u_{6} = lm$ denote the antisymmetric pairs of indices such that we obtain a $6\times 6$ matrix representation of the possible combinations of the graphs formed by the antisymmetric pairs given by
\begin{equation}\label{eq:Fuzzymatrix}
    \mathcal{G} = 
    \left( \begin{array}{ccccccc}
    G_{ikik} & G_{ikil} & G_{ikim} & G_{ikkl} & G_{ikkm} & G_{iklm} \\
    G_{ilik} & G_{ilil} & G_{ilim} & G_{ilkl} & G_{ilkm} & G_{illm} \\
    G_{imik} & G_{imil} & G_{imim} & G_{imkl} & G_{imkm} & G_{imlm} \\
    G_{klik} & G_{klil} & G_{klim} & G_{klkl} & G_{klkm} & G_{kllm} \\
    G_{kmik} & G_{kmil} & G_{kmim} & G_{kmkl} & G_{kmkm} & G_{kmlm} \\
    G_{lmik} & G_{lmil} & G_{lmim} & G_{lmkl} & G_{lmkm} & G_{lmlm} 
    \end{array}
\right)
\end{equation}
\\
\noindent
\textbf{Theorem 1.4}. The matrix of the graphs formed by the different combinations of antisymmetric indices is symmetric due to the property of union of graphs.\\
\noindent
\textit{Proof} Consider the elements $\mathcal{G}_{12} = G(i,k,i,l) = G_{ikil}$ and $\mathcal{G}_{21} = G(i,l,i,k) = G_{ilik}$ of the matrix $\mathcal{G}$. Graphically, the graphs can be expressed as a union of two subgraphs, i.e, $G(i,k,i,l) = G(i,k) \cup G(i,l)$ and $G(i,l,i,k) = G(i,l) \cup G(i,k)$. Thus, we observe that $G(i,l,i,k) = G(i,k,i,l)$ or $G_{ikil} = G_{ilik}$ and similarly, the other graphs along the principle diagonal of the matrix $\mathcal{G}$ are equal making the matrix a symmetric one.\\

\noindent
\textbf{Theorem 1.5}. The number of independent components of the subgraphs of $G$ is given by
\\
\begin{equation}
    \frac{1}{2}P(P+1) - \frac{n!}{(n-4!4!)} = \frac{n^{2}\left(n^{2} - 1 \right)}{12}.
\end{equation}
\\
\noindent
\textit{Proof} We note that for the graphs of the form $G(v_{1},v_{2},v_{1},v_{2})$, there are ${n \choose 2}$ possible choices for $v_{1}$ and $v_{2}$, for graphs of the form $G(v_{1},v_{2},v_{1},v_{3})$, there are ${n\choose 3}$ ways to choose different $v_{1}$, $v_{2}$ and $v_{2}$ and $v_{3}$ ways to choose the index that is used twice from that which results in a total of $3{n\choose 3}$ choices, and for graphs of the form $G(v_{1},v_{2},v_{3},v_{4})$, there are $2{n\choose 4}$ choices. Thus, in total we have \\
\ \ \ \ \ \ \ \ \ \ \ \ \ \ \ \ $${n\choose 2}+3{n\choose 3}+2{n\choose 4}=\frac{n^2(n^2-1)}{12}$$
\\
\noindent
Here, we have $n = 4$ indices and $4^{2}(4^{2}-1)/12 = 20$ subgraphs each accounting for the following symmetries:\\
1. $G_{iklm} = G_{lmik}$,\\
2. $G_{iklm} = -G_{kilm} = -(-G_{kiml}) = -(-(G_{ikml}))=...$ , \&\\
3. $G_{i(klm)} = 0$. \\

\noindent
\textbf{Theorem 1.6}. The graph formed by the elements of $\mathcal{G}-\{v_{PD}\}$ is a $K_{6}$ graph, where $v_{PD} =  \{u_{1}u_{1},u_{2}u_{2}, ..., u_{6}u_{6} \}$ is the vertex points of the principle diagonal of $\mathcal{G}$.\\
\noindent
\textit{Note} The $K_{6}$ graph has a vertex set, $V = \{u_{1},u_{2},u_{3},u_{4},u_{5},u_{6}\}$  and a edge set, $E = \{e_{12},e_{13},e_{15},e_{16},e_{23},e_{24},e_{25},e_{26},e_{34},e_{35},e_{36},e_{45},e_{46},e_{56}\}$, where $e_{ij} = e_{ji}$, i.e., $e_{12} = G_{ikil} = G_{ilik} = e_{21}, ..., e_{16} = G_{iklm} = G_{lmik} = e_{61},...$ \\

\noindent
\textbf{Theorem 1.7}. The number of independent components for $r$ number of permuting indices of a graph constructed from definition $1.10$ is given by:
\\
\begin{equation}
    \frac{1}{2}P(P+1) - \frac{n!}{r!(n-r)!} = \frac{1}{8}(n-1)n[(n-1)n + 2] - \frac{n!}{r!(n-r)!},
\end{equation}
\\
and since the indices $1,2,3,...,r$ permute cyclically,
\\
\begin{equation}
    G_{0(123...r)} = \frac{1}{r!}\left[G_{0123...r} + G_{0213...r} + ... + G_{0r...321} \right]
\end{equation}
\\
\noindent
\textbf{Theorem 1.8}. The elements of the principle diagonal of the matrix $\mathcal{G}$ is a set of ${}^{(r+1)}C_{2}$ number of vertices.\\
\noindent
\textit{Proof} For $G(v_{1},v_{2},v_{3},v_{4}) = G_{iklm}$, we have four indices amongst which one is fixed (which is $i$) and two slots available for the formation of a pair. Thus, to choose a pair, i.e., two indices out of four there are ${}^{4}C_{2} = 6$ available combinations. Thus, for $(r+1)$ number of vertices there are ${}^{(r+1)}C_{2}$ available combinations and since each combination was labelled as a vertex $u_{i}$, there are ${}^{(r+1)}C_{2}$ number of vertices.\\
\noindent
\textit{Note} An important property of cyclic symmetry is that this theorem holds if and only if $(r+1)$ is a number divisible by $2$. This arises due to the condition that we require an even number of indices for formation of pairs.

\section{The Fuzzy Graphical Approach}
\label{sec: Discussion on Energy Dissipation}
\textbf{Definition 2.1}[4,5]. Let $V$ be a non empty set. A fuzzy graph is a pair of functions $G = \left(\sigma,\mu \right)$, where $\sigma$ is a fuzzy subset of $V$ and $\mu$ is a symmetric fuzzy relation on $\sigma$, i.e., $\sigma: V\rightarrow [0,1]$ and $\mu: V\times V \rightarrow [0,1]$ such that $\mu(u,v)\leq \sigma(u) \wedge \sigma(v)$, for all $u,v \in V$.\\

\noindent
\textbf{Definition 2.2}. A fuzzy graph $G$ with the vertex set $V = \{v_{1}, v_{2}, v_{3}, v_{4} \}$, where $v_{1} = i, v_{2} = k, v_{3} = l, v_{4} = m$, is analogous to the Riemann tensor if it satisfies the following:\\
a. The vertex $v_{1}$ is fixed and is connected to only one other vertex, \\
b. The fixed vertex has a vertex membership, $\sigma = 1$, since the probability of finding $v_{1}$ in that vertex position is definite, \\
c. There are three vertices which span a $K_{3}$ graph amongst themselves, \\
d. The permuting vertices each have a vertex membership, $\sigma = 1/3$. since the probability of finding one of the permuting vertices at a particular vertex position of the $K_{3}$ is equal to $1/3$, \\
e. The antisymmetry of the index pairs which reflect on the type of graph is expressed via the number of cycles traversed in the $K_{3}$ graph. Let $P_{v_{2}}, P_{v_{3}}$ and $P_{v_{4}}$ represent the probabilities of finding the vertices $v_{1}, v_{2}$ and $v_{3}$ in the vertex positions of $K_{3}$. Now, the permuting combination $\{v_{2}, v_{3}, v_{4} \}$ yields:\\
1. $G_{iklm}$ if the number of cycles traversed is even in number, \& \\
2. $G_{ikml} (= -G_{iklm})$ if the number of cycles traversed is odd in number.\\
Thus, in general we can express the graph of a permuting combination with $m$ number of cycles as follows,
\\
\begin{equation}
    G(i,k,l,m) = \left(-1 \right)^{m} G_{iklm} = \left\{\begin{array}{lr}
        G_{iklm}, & \text{for \ m = even} \\
        -G_{iklm} = G_{ikml}, & \text{for\ m = odd} \end{array}\right\}.  
\end{equation}
\\

\noindent
\textbf{Definition 2.3}[11,12]. A fuzzy graph $G$ is said to be complete if $\mu(u,v) = \sigma(u) \wedge \sigma(v)$, for all $u,v \in V$.\\

\noindent
\textbf{Theorem 2.1}. The graph formed by the $\alpha$ number of permuting vertices of the fuzzy graph $G$, with one other fixed vertex, is a complete graph whose domination set consists of only one vertex which is connected to both the fixed vertex and the other $(\alpha-1)$ permuting vertices.\\
\noindent
\textit{Proof} Consider the fuzzy graph constructed based on definition $3.1$ with a permuting vertex set $V_{P} = \{v_{2}, v_{3}, v_{4} \}$ which yields $G_{iklm}$. Here, the vertex $v_{2}$ is the only one which is adjacent to both the fixed vertex $v_{1}$ and the other vertices from $V_{p} - \{v_{2}\}$. Thus, the domination set of $G$ is $D_{G} = \{v_{2}\}$.\\

\noindent
\textbf{Definition 2.4}[6,13]. An arc $(u,v)$ of a fuzzy graph $G$ is called a strong arc if $\mu(u,v) = \sigma(u) \wedge \sigma(v)$, for all $u,v \in V$.\\

\noindent
\textbf{Theorem 2.2}. The graph formed by the permuting vertices of the fuzzy graph constructed based on definition $1.10$ is a complete graph with all it's arcs being strong.\\
\noindent
\textit{Proof} For the permuting vertex set $V_{P} = \{v_{2}, v_{3}, v_{4} \}$ with vertex memberships $1/3$ each we have, in accordance to definition $2.3$, $\mu(v_{2},v_{3}) = \sigma(v_{2})\wedge \sigma(v_{3}) = 1/3 = \mu(v_{2},v_{4}) = \mu(v_{3},v_{4})$.\\

\noindent
\textbf{Definition 2.5}. Let the \textit{Levi-Civita} graph analogue $\epsilon(v_{1},x,y) = \epsilon^{ixy}$ be a graph such that 
\\
\begin{equation}
    (-1)^{m}\epsilon(i,x,y) = 
    \left\{\begin{array}{lr}
    \sum_{n=1}^{3}P_{n},\  \text{if}\ (i,x,y)\ \text{is}\ (v_{1},x,y), (x,y,v_{1}),\ \text{or}\ (y,v_{1},x);\ \text{m\ =\ even} \\
    -\sum_{n=1}^{3}P_{n},\ \text{if}\ (i,x,y)\ \text{is}\ (y,x,v_{1}), (v_{1},y,x),\ \text{or}\ (x,v_{1},y);\ \text{m\ =\ odd} \\
    \ \ \ 0,\ \text{if}\ v_{1} = x,\ \text{or}\ x=y,\ \text{or}\ y=v_{1};\ \text{m\ =\ 0\ for\ loops} 
    \end{array}\right\}\,
\end{equation}
\\
where $(x,y)$ is one of the permuting index pairs, i.e, $\{(v_{2},v_{3}),(v_{3},v_{4}), (v_{4},v_{2}) \}$, and $P_{n}$ is the probability of finding the index in a particular position. Since, we have three indices and the probability of finding an index in a particular position is equally likely, $P_{1}=P_{2}=P_{3}=1/3$,
\\
\\
\begin{equation}
    (-1)^{m}\epsilon(i,x,y) = 
    \left\{\begin{array}{lr}
    +1,\  \text{if}\ (i,x,y)\ \text{is}\ (v_{1},x,y), (x,y,v_{1}),\ \text{or}\ (y,v_{1},x);\ \text{m\ =\ even} \\
    -1,\ \text{if}\ (i,x,y)\ \text{is}\ (y,x,v_{1}), (v_{1},y,x),\ \text{or}\ (x,v_{1},y);\ \text{m\ =\ odd} \\
    \ \ \ 0,\ \text{if}\ v_{1} = x,\ \text{or}\ x=y,\ \text{or}\ y=v_{1};\ \text{m\ =\ 0\ for\ loops} 
    \end{array}\right\}\ .
\end{equation}
\\

\noindent
\textbf{Theorem 2.3}. The union of $\epsilon(v_{1},x,y)$ and the complete graph formed of permuting indices yields a loop at the common vertex, i.e., the vertex both adjacent to the fixed vertex and the other permuting ones. When the fuzzy graph $G = (\sigma,\mu)$ is expressed as the union of the complete fuzzy graph (formed of the permuting indices) and $\epsilon(v_{1},x,x)$, where $x$ is the common vertex, the vertex membership of the common vertex is reduced by
\\
\begin{equation}
    \sigma'(v_{2}) = \frac{\sigma(v_{2})}{\alpha},
\end{equation}
\\
where $\alpha$ is the number of permuting vertices. \\
\noindent
\textit{Proof} Consider the fuzzy graph $G_{i,k,l,m}$, where the vertex $v_{2}$ is the commonly adjacent to both the fixed and the permuting edges. Now,  we can define $\epsilon(v_{1},x,x) = \epsilon(i,k,k)$ and since the graph is to be expressed as the union of $\epsilon(i,k,k)$ and $G(i,k,l,m)$, we have $\epsilon(i,k,k) \cup G_{k,l,m} =  G_{i,k,l,m}$ in which $\sigma(v_{3}) = \sigma(v_{4}) = 1/3$, $\sigma(v_{1}) = 1$, and $\sigma'(v_{2}) = \frac{\sigma(v_{2})}{3} = 1/9$. \\

\section{On Route to a Fuzzy Petrov-Penrose Classification}
Generally, gravitational fields are classified in accordance to the \textit{Petrov-Penrose classification} of their corresponding Weyl tensor. This is an algebraic classification based on the idea that the curvature tensor can be thought of as a $6\times 6$ matrix and the reduction of these matrix naturally results in general categories of curvature tensors. In this section, we present the relations between the fuzzy graphical analog of the Riemann tensor and the Riemann tensor and more specifically, study the similarities of the properties between the matrices as defined in equations \ref{eq:Fuzzymatrix} and \ref{eqn:Matrix}. From the symmetries of the \textit{Riemann curvature tensor}, we can write it as a $R^{\alpha \beta}_{\gamma \delta}$ and associate an index $\mathcal{I} = 1,2,...,6$ with each pair $01, 02, 03, 23, 31, 12$ of the independent values that $\alpha \beta$ and $\gamma \delta$ can take. The curvature tensor can be expressed as a $6 \times 6$, $\mathcal{M}^{\mathcal{I}}_{\mathcal{K}}$ matrix as given below
\\
\begin{equation} \label{eqn:Matrix}
    \left( \begin{array}{ccccccc}
R^{01}_{01}\ & R^{01}_{02}\ & R^{01}_{03} & | &  R^{01}_{23}\ & R^{01}_{31}\ & R^{01}_{12} \\

R^{02}_{01}\ & R^{02}_{02}\ & R^{02}_{03} & | &  R^{02}_{23}\ & R^{02}_{31}\ & R^{02}_{12} \\

R^{03}_{01}\ & R^{03}_{02}\ & R^{03}_{03} & | &  R^{03}_{23}\ & R^{03}_{31}\ & R^{03}_{12} \\

----\ & ----\ & ---- &----|---- & ----\ & ----\ & ----\\

R^{23}_{01}\ & R^{23}_{02}\ & R^{23}_{03} & | &  R^{23}_{23}\ & R^{23}_{31}\ & R^{23}_{12} \\

R^{31}_{01}\ & R^{31}_{02}\ & R^{31}_{03} & | &  R^{31}_{23}\ & R^{31}_{31}\ & R^{31}_{12} \\

R^{12}_{01}\ & R^{12}_{02}\ & R^{12}_{03} & | &  R^{12}_{23}\ & R^{12}_{31}\ & R^{12}_{12} \end{array}
\right)
\end{equation}
\\
The matrix $\mathcal{M}^{\mathcal{I}}_{\mathcal{K}}$ can alternatively be expressed as follows
\\
\begin{equation}
    \mathcal{M}^{\mathcal{I}}_{\mathcal{K}} =  \left( \begin{array}{cc}
    \mathcal{A} & \mathcal{B}   \\
    -\mathcal{B}^{T}  & \mathcal{C} \end{array}
    \right),
\end{equation}
\\
where $\mathcal{A}, \mathcal{B}, \mathcal{C}$, are $3\times3$ matrices. Notice that the matrix $\mathcal{B}$ is null. This can be shown by first lowering the index and making use of the property of the \textit{Levi-Civita} symbol as follows
\\
\begin{equation}
    Tr\ \mathcal{B} = R^{01}_{23} + R^{02}_{31} + R^{03}_{12} = \epsilon^{011}R_{123} + \epsilon^{022}R_{231} + \epsilon^{033}R_{312}.
\end{equation}
\\
Comparing the matrix of equation \ref{eqn:Matrix} to that of equation \ref{eq:Fuzzymatrix}, we can prove that the trace of matrix $\mathcal{B}$ is null in the graph theoretical case by expressing the fuzzy graph as the union of the complete fuzzy graph and $\epsilon\left(v_{1},x,x \right)$ as follows,
\\
\begin{equation}
\begin{array}{lr}
    Tr\ \mathcal{B} = G(i,k,k,l)+G(i,l,k,m)+G(i,m,l,m)\\
    \ \ \ \ \ \ \ \ =\epsilon(i,k,k)\cup G(k,m,l)+\epsilon(i,l,l)\cup G(l,k,m) + \epsilon(i,m,m)\cup G(m,l,m)\\
    \ \ \ \ \ \ \  \ = 0.
\end{array}    
\end{equation}
\\
We also observe that the matrices $\mathcal{A}$ and $\mathcal{C}$ are equal to their transposes, i.e., $\mathcal{A} = \mathcal{A}^{T}$ and $\mathcal{C} = \mathcal{C}^{T}$. The structure of the matrix represented in equation \ref{eqn:Matrix} is based on separating the components of the Riemann curvature tensor into three distinct sets, $R_{0 \alpha 0 \alpha}$, $R_{0 \beta \gamma \delta}$, and $R_{\gamma \delta \mu \nu}$. Observe that the first set is a $3\times 3$ matrix in the indices $\alpha$ and $\beta$ and as for the other two, they are to be fixed by removal of antisymmetry that they possess. Thus, we introduce the following $3\times 3$ matrices along with their fuzzy counterpart
\\
\begin{equation}
    \begin{array}{lr}
    \Psi_{\alpha \beta} = R_{0 \alpha 0 \beta} = G(i,\alpha,i,\beta), \\ 
    \Sigma_{\alpha \beta} = \frac{1}{2}\epsilon_{\alpha \gamma \delta}R^{\gamma \delta}_{0 \beta} = \frac{1}{2}\epsilon(\alpha, v_{1},v_{2})\cup G(v_{1}, v_{2},i, \beta) \\
    \Lambda_{\alpha \beta} = \frac{1}{4}\epsilon_{\alpha \gamma 
    \delta}\epsilon_{\beta \mu \nu}R^{\gamma \delta \mu \nu} = \frac{1}{4}\epsilon(\alpha, v_{1},v_{2})\cup \epsilon(\beta,v_{3},v_{4})\cup G(v_{1},v_{2},v_{3},v_{4}),
\end{array}
\end{equation}
\\ \index{Ricci flat}
where $\epsilon_{abc}$ is a three-dimensional Levi-Civita tensor. These matrices yield the following relations under the \textit{Ricci flatness} condition, $R_{XY} = 0$
\\
\begin{equation}
    \Psi_{\alpha \alpha} = 0, \ \Sigma_{\alpha \beta} = \Sigma_{\beta \alpha}, \ \Psi_{\alpha \beta} = -\Lambda_{\alpha \beta}. 
\end{equation}
\\
According to the definitions given above we have the matrix $\Psi_{\alpha \beta}$ to have the following form 
\\
\begin{equation}
    \begin{array}{lr}
    \Psi_{11} = R_{0101}, \ \Psi_{12} = R_{0102}, \ \Psi_{0103} = R_{0103}, \ ...\\\\
    
    \Longrightarrow \Psi_{\alpha \beta} = \left( \begin{array}{ccc}
    R_{0101} & R_{0101} & R_{0103} \\
    R_{0201} & R_{0202} & R_{0203} \\
    R_{0301} & R_{0302} & R_{0303} \end{array}
    \right) = \left( \begin{array}{ccc}
    G_{ikik} & G_{ikil} & G_{ikim} \\
    G_{ilik} & G_{ilil} & G_{ilim} \\
    G_{imik} & G_{imil} & G_{imim} \end{array}
    \right)\\
    \ \ \ \ \ \ \ \ \ \ \ \ \  = \left( \begin{array}{ccc}
    \Psi_{11} & \Psi_{12} & \Psi_{13} \\
    \Psi_{21} & \Psi_{22} & \Psi_{23} \\
    \Psi_{31} & \Psi_{32} & \Psi_{33} \end{array}
    \right)
    \end{array}
\end{equation}
\\
Comparing this matrix to the $6\times 6$ form obtained previously, we find that $\Psi_{\alpha \beta}$ is comprised of the components of first quarter of the matrix (after lowering their index). Thus, $\Psi_{\alpha \beta} = \mathcal{A}$. Now, to the matrix $\Sigma_{\alpha \beta}$. Observe that in the components of the $\Sigma_{\alpha \beta}$
\\
\begin{equation}
    \Sigma_{11} = \frac{1}{2}\epsilon_{123}R_{01}^{23}, \ \Sigma_{12} = \frac{1}{2}\epsilon_{123}R_{02}^{23}, \ \Sigma_{13} = \frac{1}{2}\epsilon_{123}R_{03}^{23}, \ ...,
\end{equation}
\\
the factor $1/2$ is removed by the symmetry of the matrix, i.e., since $\Sigma_{\alpha \beta} = \Sigma_{\beta \alpha}$, $\Sigma_{12} = \Sigma_{21}, \ ...$, and hence
\\
\begin{equation}
    \Sigma_{(12)} = 2\Sigma_{12} = \Sigma_{12} + \Sigma_{21} = \underbrace{\epsilon_{123}}_{=1}R_{02}^{23} = R_{02}^{23}. 
\end{equation}
\\
Similarly, we can calculate the other components to obtain the following matrix
\\
\begin{equation}
    \Sigma_{\alpha \beta} = \left( \begin{array}{ccc}
    R_{01}^{23} & R_{02}^{23} & R_{03}^{23} \\
    R_{01}^{31} & R_{02}^{31} & R_{03}^{31} \\
    R_{01}^{12} & R_{02}^{12} & R_{03}^{12} \end{array}
    \right) = \left( \begin{array}{ccc}
    G_{lmik} & G_{lmil} & G_{lmim} \\
    G_{mkik} & G_{mkil} & G_{mkim} \\
    G_{klik} & G_{klil} & R_{klim} \end{array}
    \right) = \left( \begin{array}{ccc}
    \Sigma_{11} & \Sigma_{12} & \Sigma_{13} \\
    \Sigma_{21} & \Sigma_{22} & \Sigma_{23} \\
    \Sigma_{31} & \Sigma_{32} & \Sigma_{33} \end{array}
    \right).
\end{equation}
\\ \index{block symmetry}
Comparing this matrix to the $6\times 6$ form obtained previously, we find that $\Sigma_{\alpha \beta}$ is comprised of the components of third quarter of the matrix (, i.e. the first half of the second row). Thus, $\Sigma_{\alpha \beta} = -\mathcal{B}^{T}$. In matrix $\Lambda_{\alpha \beta}$, notice that there is symmetry in the indices and also among matrix components due to the \textit{block symmetry} of the curvature tensor. the following are the components of the matrix $\Lambda_{\alpha \beta}$
\\
\begin{equation}
\begin{array}{lr}
    \Lambda_{11} = \frac{1}{4}\epsilon_{123}\epsilon_{123}R^{1323}, \ \Lambda_{12} = \frac{1}{4}\epsilon_{123}\epsilon_{231}R^{2331}, \ \Lambda_{13} = \frac{1}{4}\epsilon_{312}\epsilon_{312}R^{1212}, \\\\
    
    \Lambda_{21} = \frac{1}{4}\epsilon_{213}\epsilon_{123}R^{1323}, \ ...
\end{array}
\end{equation}
\\
We know that $\Lambda_{\alpha \beta}$ is a symmetric matrix thus, components such as $\Lambda_{12} = \Lambda_{21} \Longrightarrow \Lambda_{(12)} = 2 \Lambda_{12}$, and this eliminates the factor $(1/2)$. Now, to account for the remaining $(1/2)$, consider the matrix components $a_{12} = \Sigma_{12}$ and $a_{21} = \Sigma_{21}$ (using index $a$ to avoid confusion), in which there exists a block symmetry\footnote{$R^{\alpha \beta \gamma \delta} = R^{\delta \gamma \alpha \beta} = R^{\beta \alpha \gamma \delta}$} between Riemann curvature tensor components, $R^{2331} = R^{1323}$. This implies that
\\
\begin{equation}
\begin{array}{lr}
    a_{12} = \Sigma_{12} = \frac{1}{2}\epsilon_{123}\epsilon_{231}R^{2331} = \frac{1}{2}\epsilon_{213}\epsilon_{123}R^{1323} = \Sigma_{21} = a_{21}\\\\
    
    \Longrightarrow 2 a_{(12)} = a_{12} + a_{21} = \underbrace{\epsilon_{123}}_{=1}\underbrace{\epsilon_{231}}_{=1}R^{2331} = R^{2331}.
\end{array}
\end{equation}
\\
Similarly, we can calculate the other components to obtain the following matrix
\\
\begin{equation}
    \Lambda_{\alpha \beta} = \left( \begin{array}{ccc}
    R^{2323} & R^{2331} & R^{2312} \\
    R^{3123} & R^{3131} & R^{3112} \\
    R^{1223} & R^{1231} & R^{1212} \end{array}
    \right) = \left( \begin{array}{ccc}
    G_{lmlm} & G_{lmmk} & G_{lmkl} \\
    G_{mklm} & G_{mkmk} & G_{mkkl} \\
    G_{kllm} & G_{klmk} & G_{klkl} \end{array}
    \right)= \left( \begin{array}{ccc}
    \Lambda_{11} & \Lambda_{12} & \Lambda_{13} \\
    \Lambda_{21} & \Lambda_{22} & \Lambda_{23} \\
    \Lambda_{31} & \Lambda_{32} & \Lambda_{33} \end{array}
    \right)
\end{equation}
\\
Comparing this matrix to the $6\times 6$ form obtained previously, we find that $\Lambda_{\alpha \beta}$ is comprised of the components of fourth quarter of the matrix (, i.e. the second half of the second row). Thus, $\Lambda_{\alpha \beta} = \mathcal{C}$.
\\\\
Let $\Omega_{\alpha \beta}$ be a symmetric complex tensor defined as follows
\\
\begin{equation}
\begin{array}{lr}
    \Omega_{\alpha \beta} = \frac{1}{2}\left(\Psi_{\alpha \beta} + 2i\Sigma_{\alpha \beta} - \Lambda_{\alpha \beta} \right) = \frac{1}{2}\left(\Psi_{\alpha \beta} + 2i\Sigma_{\alpha \beta} + \Psi_{\alpha \beta} \right) = \Psi_{\alpha \beta} + i\Sigma_{\alpha \beta}\\
    \Omega_{\alpha \beta} = G(i,\alpha,i,\alpha) + \frac{i}{2}\epsilon(\alpha,v_{1},v_{2})\cup G(v_{1},v_{2},i,\beta).
\end{array}
\end{equation}
\\
It is well known that classification of the Riemann curvature tensor can be reduced to a simple eigen value problem where we consider the eigen value equation $\Omega_{\alpha \beta} k_{\beta} = \lambda k_{\alpha}$, in which the complex eigenvalues $\lambda = \lambda_{R} + i\lambda_{I}$ satisfy the condition $\lambda^{(1)}+ \lambda^{(2)} + \lambda^{(3)} = 0$ since $\Omega_{\alpha}^{\alpha} = 0$. The matrix's classification is now dependent on the number of independent eigenvectors and leads to six different cases, called Petrov Types $I$, $II$, $D$, $III$, $N$, and $O$.

\section{Conclusion}
In this paper we have introduced a graph theoretic and fuzzy graph theoretic analog of the Riemann tensor and have shown how the latter satisfies the properties of the $6\times 6$ matrix of the Riemann tensor. We then use the definitions and the properties discussed in sections $2$ and $3$ to develop a fuzzy graphical analog of the Petrov-Penrose classification. We hope to study the detailed fuzzy graphical connections among the various Petrov types in future papers where we would explore connections between fuzzy graphs and exterior algebra and study fuzzy analogs of the \textit{Weyl tensor}, the \textit{Hodge star operator}, and the \textit{Kretchmann invariant}.

\end{document}